\documentclass{article}[10]

\usepackage{amssymb}
\usepackage{latexsym}
\usepackage{amsmath}
\usepackage[mathscr]{eucal}
\usepackage{epsfig}

\newtheorem{Lemma}{Lemma}
\newtheorem{Theorem}{Theorem}
\newtheorem{Remark}[Lemma]{Remark}
\newtheorem{Proposition}[Theorem]{Proposition}

\newtheorem{Assumption}[Lemma]{Assumption}

\newcommand{\CC}{\mbox{${\mathcal C}$}}

\newcommand{\GG}{\mbox{${\mathcal G}$}}
\newcommand{\FF}{\mbox{${\mathcal F}$}}

\newcommand{\HH}{\mbox{${\mathcal H}$}}
\newcommand{\WW}{\mbox{${\mathcal W}$}}

\renewcommand{\SS}{\mbox{${\mathcal S}$}}

\newcommand{\Nbold}{\mbox{${\mathbb N}$}}
\newcommand{\Rbold}{\mbox{${\mathbb R}$}}

\newcommand{\Zbold}{\mbox{${\mathbb Z}$}}
\newcommand{\Lbold}{\mbox{${\mathbb L}$}}

\newcommand{\Ebb}{\mbox{${\mathbb E}$}}
\newcommand{\Pbb}{\mbox{${\mathbb P}$}}

\newcommand{\BBB}{\mbox{${\mathscr{B}}$}}

\newcommand{\bE}{{\bf E}}
\newcommand{\bP}{{\bf P}}

\newcommand{\bQ}{{\bf Q}}

\newcommand{\ba}{{\bf a}}

\newcommand{\be}{{\bf e}}
\newcommand{\bx}{{\bf x}}
\newcommand{\by}{{\bf y}}

\newcommand{\bv}{{\bf v}}
\newcommand{\bz}{{\bf z}}

\newcommand{\bzero}{{\bf 0}}
\newcommand{\bone}{{\bf 1}}

\newcommand{\Var}{{\rm Var}\ }
\newcommand{\eps}{\varepsilon}
\newcommand{\sfrac}[2]{{\textstyle\frac{#1}{#2}}}

\newcommand{\qed}{\hfill \ \ \rule{1ex}{1ex}}
\newcommand{\proof}{\noindent\emph{Proof :}\ \ \ \ }

\newcommand{\Pbar}{{\overline \Pbb}}
\newcommand{\Ebar}{{\overline \Ebb}}

\newcommand{\Ybar}{{\overline Y}}

\title{{\bf Random Walk in Dynamic Markovian Random Environment}}

\author{{\bf Antar Bandyopadhyay}\footnote{email:
antar@math.chalmers.se. This work was begun while 
A.B. was a post-doctoral fellow at the IMA, University of Minnesota.} 
\vspace{0.1in}\\
        Department of Mathematics\\
        Chalmers University of Technology\\
        SE-412 96, G\"{o}teborg, SWEDEN \vspace{0.1in}
        \and
        {\bf Ofer Zeitouni}\footnote{email: zeitouni@math.umn.edu.
Partially supported
by NSF grants DMS-0302230 and DMS-0503775} \vspace{0.1in}\\
        Department of Mathematics \\
        University of Minnesota \\
        Minneapolis 55455, USA \\
        and\\
Technion, Haifa, Israel}

\date{March 27, 2006}
\begin{document}

\maketitle
\bibliographystyle{plain} 

\begin{abstract}
We consider a model, 
introduced by Boldrighini, Minlos and 
Pellegrinotti \cite{BMP97, BMP00},
of
random walks in \emph{dynamical} random environments
on the integer lattice $\Zbold^d$ with $d \geq 1$. In this model,
the environment changes over time
in a Markovian manner, independently across sites, while the
walker uses the environment at its current location in order
to make the next transition. 
In contrast with the cluster expansions approach of
\cite{BMP00}, we follow 
a probabilistic argument
based on regeneration times.
We prove
an annealed SLLN and invariance principle for any dimension,
and provide a quenched invariance principle for dimension $d > 7$, 
providing for $d>7$ an
alternative to the analytical approach
of \cite{BMP00}, with the added benefit 
that it is valid under weaker assumptions.
The quenched results use, in addition to the regeneration
times already mentioned, a
technique introduced by   Bolthausen and Sznitman 
\cite{BS}.
\end{abstract}



\newpage

\section{Introduction}
\label{Sec:Intro}

In recent years there has been a great deal of study on 
\emph{random walks in random environments} (RWRE) on 
$\Zbold^d$, $d \geq 1$, where
first the \emph{environment} is chosen at random and kept fixed 
throughout the time evolution, and then a walker moves randomly in such a way
that
given the environment, its position forms
a time homogeneous Markov chain whose transition
probabilities depend on the environment. 
Even though a lot is known about this model,
there are many challenging problems left open, 
see \cite{Zei04} and \cite{sznitmanLN}
for  surveys.  

In the current paper, we consider random walks in {\it dynamical}
random environments,
where along with the walker the environment changes 
over time. Such model was studied in an abstract setting
by Kifer \cite{Ki96}. More relevant to us, the model we consider was
first introduced by 
Boldrighini, Minlos and Pellegrinotti \cite{BMP97}, where they 
consider the case where the environment changes over time
in an \emph{i.i.d.} fashion. In \cite{BMP97} and in
\cite{BMP03}, they proved under certain assumptions (omitted here)
that for almost every realization of the dynamic environment,
the position of the random walk, properly centered and scaled,
satisfies a \emph{central limit theorem},
with covariance that 
does not depend on the particular realization
of environment. 
Further, in
\cite{BMP00}, they obtain similar results for $d\geq 3$
for certain environment evolving over time as a Markov chain, 
independently at 
each site. 
The proofs in these papers are based on
 \emph{cluster expansion}, and involve a heavy analytic machinery.  
In the i.i.d. case,
a somewhat simpler proof 
can be found in Stannat \cite{Sta04}. 

Our goal in this paper
is to describe a probabilistic treatment 
of this model, which is arguably simpler than those that
have appeared in the literature. 
We recover most of the results of \cite{BMP00}, at least when
$d>7$, under weaker and more natural hypotheses.
Our approach is based on the introduction of appropriate
``regeneration times'', borrowing this concept, if not
the details of the construction, from the RWRE literature.
Our approach to quenched results is based on a technique introduced
in \cite{BS}.

The paper is divided as follows. The next section describes precisely 
the model and states our main results: a strong law of large numbers,
an annealed CLT (in any dimension), and a quenched CLT (for $d>7$).
Section \ref{sec-regeneration} constructs the regeneration
times alluded to above, and derives their basic properties.
Section \ref{Sec:Proofs} provides the proofs of our main results.
Finally, we present in Section \ref{Sec:Remarks} several remarks 
and open problems.
\section{Description of the model and main results}
\label{Sec:Main-Results} 

In what follows we will consider $\Zbold^d$, for a fixed $d \geq 1$ with
nearest neighbor graph structure as our basic underlying graph. 
Let 
$N_{\bzero} := \left\{\pm \be_i\right\}_{i=1}^d \cup \left\{ \bzero \right\}$
be the neighbors of the origin $\bzero$ (including itself), where
$\be_i$ is the unit vector in the $i^{\mbox{th}}$ coordinate direction. 
We use 
$N_{\bx} := \bx + N_{\bzero}$ to denote the collection
of neighbors of an  
$\bx \in \Zbold^d$. 

Let $\SS := \SS_{\bzero}$ 
be a collection of probability
measures  on the $\left(2d+1\right)$ elements of $N_{\bzero}$. 
To simplify the presentation and avoid various measurability issues,
we assume that $\SS$ is a Polish space (including the possibilities
that $\SS$ is finite or countably infinite).
For each $\bx \in \Zbold^d$,  $\SS_{\bx}$ denotes  a copy of $\SS$,
with all elements of $\SS$ shifted to have support on
$N_{\bx}$. Formally, an element
$\omega(\bx,\cdot)$ of $\SS_{\bx}$, is a probability measure
satisfying
\[
\begin{array}{rcl}
\omega\left(\bx, \by\right) \geq 0 \,\,\,\, \forall \,\,\, 
\bx, \by \in \Zbold^d
& 
\mbox{and}
& 
\mathop{\sum}\limits_{\by \in N_{\bx}} \omega\left(\bx, \by\right) = 1
\end{array}
\]


Let ${\cal B}_{\SS}$
denote the Borel $\sigma$-field on $\SS$. We let
$K \colon \SS \times {\cal B}_{\SS} \rightarrow [0,1]$ denote
a Markov transition
probability on $\SS$ with a unique stationary distribution $\pi$. Let 
$P^{\pi}$ be the probability distribution on the standard product space
$\SS^{\Nbold_0}$, where $\Nbold_0 := \left\{ 0 \right\} \cup \Nbold$, 
giving the $\pi$-stationary Markov $K$-chain on $\SS$. 
For each $\bx \in \Zbold^d$, let $K_{\bx}$, 
$\pi_{\bx}$ and $P^{\pi_{\bx}}$ be just the copies of respective
quantities when we replace $\SS$ by $\SS_{\bx}$. Let 
\begin{equation}
\begin{array}{cc}
\Omega := \mathop{\prod}\limits_{\bx \in \Zbold^d} \SS_{\bx}^{\Nbold_0},
& 
\bP^{\pi} := \mathop{\bigotimes}\limits_{\bx \in \Zbold^d} P^{\pi_{\bx}},
\end{array}
\label{Equ:Omega-P}
\end{equation}
where the measure is defined on $\FF$, 
the standard product $\sigma$-algebra on $\Omega$. 
An element $\omega \in \Omega$ will be written as
$\left\{ \left( \omega_t\left(\bx, \cdot\right) \right)_{t \geq 0} 
\,\Big\vert\, \bx \in \Zbold^d \,\right\}$.
We note
that under $\bP^{\pi}$, 
the canonical variable $\omega \in \Omega$ follows
a distribution such that 
\begin{itemize}
\item For each $\bx \in \Zbold^d$, 
      $\left(\omega_t\left(\bx, \cdot\right)\right)_{t \geq 0}$ is a 
      stationary Markov chain with transition kernel $K$.\footnote{ 
      By this we mean that 
      $\left(\omega_t\left(\bx, \cdot\right)\right)_{t \geq 0}$ is a 
a stationary Markov chain
      on $\SS_{\bx}$ with transition probabilities 
$K_{\bx}$ and stationary law
      $\pi_{\bx}$. 
To alleviate notation, we identify in the sequel
$\SS_{\bx}$  with $\SS$, $\pi_{\bx}$ with $\pi$, and $K_{\bx}$ with
$K$ whenever no confusion occurs.}
\item The chains $\left(\omega_t\left(\bx, \cdot\right)\right)_{t \geq 0}$
      are i.i.d. as $\bx$ varies over $\Zbold^d$. 
\end{itemize}

We now turn to define a  random walk $(X_t)_{t\geq 0}$.
Given an environment 
$\omega \in \Omega$, $(X_t)_{t\geq 0}$
is a time inhomogeneous Markov chain
taking values in $\Zbold^d$ 
with transition probabilities
\begin{equation}
\bP_{\omega} \left( X_{t+1} = \by \,\Big\vert\, X_t = \bx \right)
= \omega_t\left(\bx, \by\right). 
\label{Equ:RW-Trans-Prob}
\end{equation}
For each $\omega \in \Omega$, we denote by
$\bP_{\omega}^{\bx}$ the law induced by $(X_t)_{t\geq 0}$ on
$\left( \left(\Zbold^d\right)^{\Nbold_0}, \GG \right)$,
where $\GG$ is the $\sigma$-algebra generated by the cylinder sets, 
such that
\begin{equation}
\bP_{\omega}^{\bx}\left( X_0 = \bx \right) = 1. 
\label{Equ:RW-Intitial-Prob}
\end{equation}
Naturally $\bP_{\omega}^{\bx}$ is called the \emph{quenched law}
of the random walk $\left\{X_t\right\}_{t \geq 0}$, starting at $\bx$. 

We note that for every $G \in \GG$, the function
\[
\omega \mapsto \bP_{\omega}^{\bx}\left(G\right)
\]
is $\FF$-measurable. Hence, we may define the measure
$\Pbb^{\bx}$ on 
$\left( \Omega \times \left(\Zbold^d\right)^{\Nbold_0},
\FF \otimes \GG \right)$ from the relation
\[
\Pbb^{\bx}\left( F \times G \right)
= 
\int_F \! \bP_{\omega}^{\bx}\left(G\right) \bP^{\pi}\left(d\omega\right), 
\,\,\,\, \forall \,\,\, F \in \FF, \, G \in \GG. 
\]
With a slight abuse of notation, we also denote the marginal of 
$\Pbb^{\bx}$ on $\left(\Zbold^d\right)^{\Nbold_0}$
by $\Pbb^{\bx}$, whenever no confusion occurs. This probability distribution
is called
the \emph{annealed law} of the random walk $\left\{X_t\right\}_{t \geq 0}$, 
starting at $\bx$. 
Note that under
$\Pbb^{\bx}$, the random walk 
$\left\{ X_t \right\}_{t \geq 0}$ 
is not, in general, a Markov chain.
However, when $K(s,\cdot)=K(s_0,\cdot)$ for some $s_0\in \SS$,
i.e. the environment is i.i.d.
in time, 
then of course
$\left\{X_t\right\}_{t \geq 0}$ is actually a random walk on $\Zbold^d$
under $\Pbb^{\bx}$, with deterministic increment distribution given by
the mean of $\pi$. 

Throughout this paper we will assume the followings hold,
\begin{Assumption}
\label{Ass-1}

\begin{itemize}
\item[(A1)] There exists $0 < \kappa \leq 1$ such that 
            \begin{equation}
            K\left(w, A\right) \geq \kappa \, \pi\left(A\right),
            \,\,\,\, \forall \,\,\, w \in \SS, A\in {\cal B}_{\SS}.
            \label{Equ:Assumption-1}
            \end{equation}   
\item[(A2)] There exist $0 < \eps < 1$ and 
            a fixed translation invariant
Markov kernel with only nearest neighbor transition
            $q \colon \Zbold^d \times \Zbold^d \rightarrow [0,1]$ 
            with the property that $q(\bx,\by)=q(\by-\bx)$ and
            \begin{equation}
            \Big\vert \sum_{\by \in \Zbold^d:|y|=1} 
q\left(0,\by\right) \, 
            e^{\imath \, \ell \cdot \by} \Big\vert < 1, 
            \,\,\,\, \forall \,\,\,
            \ell \in \Zbold^d \setminus \left\{\bzero\right\},
            \label{Equ:q-non-trivial}
            \end{equation}  
            such that 
            \begin{equation}
            \bP^{\pi}\left( \omega_t\left( \bx, \by \right) \geq 
            \eps \, q\left(\bx, \by\right) \right) = 1, 
            \,\,\,\, \forall \,\,\, \bx, \by \in \Zbold^d, \,\,\, t \geq 0.
            \label{Equ:Uniform-Elliptic}
            \end{equation}
\item[(A3)] $\kappa + \eps^2 > 1$. 
\end{itemize}
\end{Assumption}
\begin{Remark}
Some comments regarding Assumption \ref{Ass-1} are in order.
\begin{enumerate}

\item
Condition (A1) 
provides a uniform ``fast mixing'' rate for the environment
chains. If, as in \cite{BMP00},
$\SS$ is finite and $K$ is irreducible and aperiodic, then while
(A1) may fail, it does hold if $K$ is replaced by $K^r$ 
for some fixed $r\geq 1$. A slight modification of our arguments applies
to that
case, too. 

\item
Condition (A2)  is the
same as that made in \cite{BMP97, BMP00} and later on in \cite{Sta04}. 
This condition essentially means that
the random environment  has 
a ``deterministic'' part $q$, which is non-degenerate. We thus refer
to condition (A2) as an
\emph{ellipticity} condition. 

\indent
We note that the assumption that $q$ is translation 
invariant does not present a 
loss of generality.
This is because under $\bP^{\pi}$ the
environment chains $\left(\omega_t\left(\bx, 
\cdot\right)\right)_{t \geq 0}$
are assumed i.i.d. as $\bx$ varies. 

\indent
Finally note that we could allow 
$\eps=1$ in (\ref{Equ:Uniform-Elliptic}). In that case,
the environment is given
by the deterministic kernel $q$. 
This makes the walk a classical Markov
chain which is well studied. This is the
reason why we work with the restriction  $ \eps < 1$.
 
\item Condition (A3) is
technical but absolutely crucial for our argument. It implies that
there is a trade off between 
the environment Markov chain being fast mixing
($\kappa$  close
to $1$) and the fluctuation in the environment being ``small'' ($\eps$ close
to $1$). The later is a condition assumed in \cite{BMP97, BMP00} also. 

\end{enumerate}
\end{Remark}
We now formulate our main results as follows. 
\begin{Theorem}[Annealed SLLN]
\label{Thm:SLLN}
Let Assumption \ref{Ass-1} hold.
Then
there exists a deterministic $\bv \in \Rbold^d$ such that 
\begin{equation}
\frac{X_n}{n} \mathop{\longrightarrow}\limits_{n\to\infty} \bv \,\,\,\,
\mbox{a.s.\ \ } \left[\Pbb^{\bzero}\right]
\label{Equ:SLLN}
\end{equation}
\end{Theorem}
\begin{Remark}
A (non-explicit) formula for the limit velocity $\bv$ is given below, see
(\ref{Equ:Final-Limit}). A consequence of that
formula is that whenever
the transition probabilities $K(s,\cdot)$ are invariant under lattice 
isometries (i.e., $\int_{\SS}K(s,ds') s'(e)=
 \int_{\SS}K(s,ds') s'(Te)$ for any lattice isometry $T$ and
any $s\in \SS$), then
$\bv=0$.
\end{Remark}%
%
We also prove
\begin{Theorem}[Annealed Invariance Principle]
\label{Thm:CLT}
Let Assumption \ref{Ass-1} hold.
Then,
there exists a $\left( d \times d \right)$ 
strictly  positive definite matrix $\Sigma$,
such that under $\Pbb^{\bzero}$,
\begin{equation}
\left(\frac{X_{\lfloor nt\rfloor } - nt \, \bv}{\sqrt{n}}\right)_{t\geq 0} 
\mathop{\longrightarrow}\limits^{d}_{n\to\infty} 
\mbox{\emph{BM}}_d\left(\Sigma\right)\,,
\label{Equ:CLT}
\end{equation}
where
$
\mbox{\emph{BM}}_d\left(\Sigma\right)$
denotes  a $d$-dimensional
Brownian motion with 
covariance matrix $\Sigma$,
$\mathop{\longrightarrow}\limits^{d}$ denotes weak convergence of the laws,
and $\lfloor x\rfloor$ denotes the integer part of $x$.
\end{Theorem}
\begin{Remark}
An implicit
formula for the covariance matrix $\Sigma$ 
is given in 
(\ref{eq-250805aa}).
\end{Remark} 
In order to get a quenched invariance principle, we need some further 
restrictions.  Set 
\begin{equation}
\label{eq-200805na}
\gamma=\log(1-\kappa)/\log \varepsilon\,.
\end{equation}
\begin{Theorem}[Quenched Invariance Principle]
\label{Thm:CLTquenched}
Let Assumption \ref{Ass-1} hold. Assume further
that $\gamma>6$  and also
\begin{equation}
\label{eq-200805a}
d > 1 + \frac{4+\frac{2\gamma (d-1)}{\gamma(d-3)-2(d-1)}
+ \frac{8(d-1)}{\gamma(d-3)}}{1-6/\gamma}> 7\,.
\end{equation}
Then, with the notation of Theorem \ref{Thm:CLT}, one has
that (\ref{Equ:CLT}) holds true under
$\bP^{0}_\omega$,
for $\bP^{\pi}$ almost every $\omega$.
\end{Theorem}
\begin{Remark}
\label{rem-200805}
Condition (\ref{eq-200805a}) is clearly not optimal, as the case
$\kappa=1$ demonstrates (recall that when $\kappa=1$,
the environment is i.i.d. in time, and  the quenched
CLT statement holds true in any dimension as soon as
(A2) holds true, see \cite{Sta04}). We do not know however
whether the quenched CLT holds for Markov environments in
low dimension.
\end{Remark}

\section{Construction of a ``regeneration time''}
\label{sec-regeneration}
Since for $\kappa=1$ the annealed law on the random walk is Markovian
with deterministic transition, we may and will in the sequel
consider only the case $\kappa<1$.
Our approach is based on the construction
of  a sequence of a.s. finite
(stopping) times $\left\{ \tau_n \right\}_{n \geq 1}$ 
(on an enlarged probability space),
such that, between two successions, say $\tau_n$ and $\tau_{n+1}$, 
environment chains at each location visited by the walk go through a 
time of ``regeneration'', in the sense that they have started
afresh from new states selected according to the stationary distribution
$\pi$. This  in turn provides a \emph{renewal} structure for
$\left\{ \tau_n, X_{\tau_n} \right\}_{n \geq 1}$.
We note that regeneration times have been extensively
used in the RWRE context, see \cite{sznitmanLN,Zei04}.
What makes the situation considered here particularly simple is that 
the regeneration times constructed here are actually stopping times.

To introduce the regeneration times,
we begin by constructing an extension of our
probability space, which is similar to what is done in
\cite{Zei04}. 

Let $W=\left\{0,1\right\}$ and let $\WW$ denote the $\sigma$-algebra on 
$W^{\Nbold}$ generated by cylinder sets. 
For  $\eps > 0$ as in Condition
(A2), let $\bQ_{\eps}$ be the product measure on 
$\left(W^{\Nbold}, \WW\right)$ such that the coordinate variables, 
say $\left(\epsilon_t\right)_{t \geq 1}$, are i.i.d. Bernoulli variables
with 
$\bQ_{\eps}\left(\epsilon_1 = 1\right) = \eps$.  

For any
$\left(\omega, {\mathbf \epsilon}\right)\in \Omega\times W^{\Nbold}$,
we define a probability measure
$\Pbar_{\omega, \eps}^{\bzero}$ on 
$\left(\Zbold^d\right)^{\Nbold_0}$ such that 
$\left\{ X_t \right\}_{t \geq 0}$ is a Markov chain with state space 
$\Zbold^d$ and transition law
\begin{equation}
\Pbar_{\omega, \eps}^{\bzero} \left( X_{t+1} = \by 
                                  \,\Big\vert\, X_t = \bx \right) 
= 
\bone_{\left[ \epsilon_{t+1} = 1\right]} \, q\left(\bx, \by\right) + 
\frac{\bone_{\left[ \epsilon_{t+1} = 0 \right]}}{1 - \eps} \,
\left[ \omega_t\left( \bx, \by \right) - \eps q\left(\bx, \by\right) \right], 
\label{Equ:Epsilon-Coin-Defi-1}
\end{equation}
where 
$\by \in \left\{\bx\right\} \cup \left\{ \bx \pm \be_i \right\}_{i=1}^d$, 
and also 
\begin{equation}
\Pbar_{\omega, \eps}^{\bzero}\left( X_0 = \bzero \right) = 1.
\label{Equ:Epsilon-Coin-Defi-2}
\end{equation} 
Finally we define a measure $\Pbar^{\bzero}$ on
$\left( \Omega \times W^{\Nbold} \times \left(\Zbold^d\right)^{\Nbold_0}, 
\FF \otimes \WW \otimes \GG \right)$, such that the coordinate variables has 
the following distribution
\[
\left( \omega, {\mathbf \epsilon} \right) \sim \bP^{\pi} \otimes \bQ_{\eps}
\mbox{\ \ and\ \ } \left\{ X_t \right\}_{t \geq 0} \sim 
\Pbar_{\omega, {\mathbf \epsilon}}^{\bzero} \mbox{\ \ given\ \ } 
\left( \omega, {\mathbf \epsilon} \right). 
\]
It is immediate to check that under $\Pbar^{\bzero}$, 
the marginal distribution of the chain $\left\{X_t\right\}_{t \geq 0}$ is 
$\Pbb^{\bzero}$ and also the conditional distribution 
of $\left\{ X_t \right\}_{t \geq 0}$ given $\omega$ is
$\bP_{\omega}^{\bzero}$. 
We will call the $\epsilon$-variables the $\epsilon$-coin tosses. 
Heuristically, under $\Pbar^{\bzero}$, the evolution of the random walk
at time $t$ is done by first tossing
the $\eps$-coin $\eps_t$, and depending on
the outcome,  taking a step either according to how the
environment dictates it (if $\eps_t=0$)
or taking
a step according 
to the fixed transition kernel $q$ (if $\eps_t=1$).

One more extension of the probability space is needed before we can 
define the ``regeneration time'' and it is for the Markov evolution of the
environment at each sites. Notice under Condition (A1), a step in the
$K$-chain on $\SS$ can be taken as follows. First, toss a coin independently 
with probability $\kappa$ of turning head, if it turns out head, then select 
a state according to the stationary distribution $\pi$ independently; 
otherwise if the coin lands in tail, then take a step according to the
Markov transition kernel 
\[ 
\tilde{K} \left( w, A \right) = 
\frac{K\left(w,A\right) - \kappa \pi\left(A\right)}{1 - \kappa}\,,
\quad w\in \SS, A\in {\cal B}_{\SS}.
\]
A rigorous construction is as follows. Extend
the measurable space 
$\left( \Omega, \FF \right)$  to accommodate i.i.d. 
$\mbox{Bernoulli}\left(\kappa\right)$ variables
$\left\{ \left( \alpha_t\left(\bx\right) \right)_{t \geq 1} \,\Big\vert\,
\bx \in \Zbold^d \,\right\}$, such that under $\bP^{\pi}$ the Markov
evolution of the environment 
$\left(\omega_t\left(\bx, \cdot\right)\right)_{t \geq 0}$ at a site 
$\bx \in \Zbold^d$ is obtained by the above description using the
$\kappa$-coin tosses $\left(\alpha_t\left(\bx\right)\right)_{t \geq 1}$. 
We note that in this construction 
\begin{equation}
\left\{ \left( \alpha_t\left(\bx\right) \right)_{t \geq s+1} 
\,\Big\vert\, \bx \in \Zbold^d \,\right\}
\mbox{\ is independent of\ } 
\left\{ \left(\omega_t\left(\bx, \cdot\right)\right)_{0 \leq t \leq s} ; 
         \left(X_t\right)_{0 \leq t \leq s} \right\},
\label{Equ:alpha-independence}
\end{equation}
for every $s \geq 0$. 

For $t \geq 0$ and $\bx \in \Zbold^d$,  define 
\begin{equation}
I_t\left(\bx\right) := 
\sum_{s=0}^t {\mathbf 1}\left( X_s = \bx, \epsilon_s = 0 \right),
\label{Equ:I-Defi}
\end{equation}
which is the number of ``proper'' visits to the site $\bx$ by the chain up to
time $t$ (``proper'' visits are those visits in which 
the walker ``learns'' about the environment, that is those in which
the next move depends on the random environment rather than on the
auxiliary $\eps$-coin). If $I_t\left(\bx\right) > 0$ let 
\begin{equation}
\gamma_t\left(\bx\right) :=
\sup \left\{ s \leq t \,\Big\vert\, X_s = \bx , \epsilon_s = 0 \, \right\}, 
\label{Equ:gamma-Defi}
\end{equation}
be the time of last ``proper'' visit to $\bx$ before time $t$, and 
\begin{equation}
\eta_t\left(\bx\right) := 
\inf \left\{ s \geq 0 \,\Big\vert\, 
\alpha_{\gamma_t\left(\bx\right)+s}\left(\bx\right) = 1
\,\right\},
\label{Equ:eta-Defi}
\end{equation}
be the first time after $\gamma_t\left(\bx\right)$ when the environment chain
at site $\bx$ takes a step according to the stationary measure $\pi$. 
For completeness we will take
$\gamma_t\left(\bx\right) = \eta_t\left(\bx\right) = 0$ if
$I_t\left(\bx\right) = 0$. 
Finally we define 
\begin{equation}
\tau_1 :=
\inf \left\{ t > 0 \,\Big\vert\, 
\gamma_t\left(\bx\right) + \eta_t\left(\bx\right) < t 
\,\,\, \forall \,\, \bx \in \Zbold^d \,\right\}. 
\label{Equ:tau-Defi}
\end{equation}
$\tau_1$ will be called a ``regeneration time'', because from time $\tau_1$
the environment chains at all sites look like they have started afresh 
from stationary distribution. 

\begin{Proposition}
\label{Prop:tau-Finite}
Let Assumption \ref{Ass-1} hold. Then,
$\Ebar^{\bzero}\left( \tau_1^2\right) < \infty\,. $
\end{Proposition}

\proof Due to the strict inequality in
Condition (A3), we can find $0 < \delta < 1$ such
that $\eps^2 > \left(1-\kappa\right)^{\delta}$. 
Define
\begin{equation}
L\left(t\right) = \lfloor 
                   - \frac{\delta \log t}{\log \eps}
                   \rfloor,
\label{Equ:Choice-of-L-t}
\end{equation}
which is an increasing sequence of integers going to $\infty$ with
$L(t) < t$ for large $t$. 

For  fixed $t \geq 1$, let $\beta_t$ be the first time
there is a run of length $L\left(t\right)$ of non-zero $\epsilon$-coin
tosses ending at it, that is
\begin{equation}
\beta_t := \inf \left\{ s \geq L\left(t\right) \,\Big\vert\,
\epsilon_s =1, \epsilon_{s-1} =1, \cdots , 
\epsilon_{s - L\left(t\right) + 1} =1 \,\right\}\,.
\label{Equ:Defi-beta-t}
\end{equation}
From the definition of $\tau_1$ we get 
\begin{eqnarray}
\label{Equ:tau-Finite-Main-Esti}
\Pbar^{\bzero}\left( \tau_1 > t \right)
& \leq &
\Pbar^{\bzero}\left( \beta_t > t \right)  \\
&    & 
\!\!\!\!\!\! + \, 
\Pbar^{\bzero}\left( \beta_t \leq t, \,\, \exists \,\,
                  \bx \in \Zbold^d
                  \mbox{\ s.t.\ }
                  \eta_{\beta_t}\left(\bx\right) \geq \beta_t
                  - \gamma_{\beta_t}\left(\bx\right) \right)
\nonumber
\end{eqnarray}
We will consider the first and second terms in the right hand side of 
(\ref{Equ:tau-Finite-Main-Esti}) separately, as follows.

Consider first the second term,  noting that for each 
$\bx \in \Zbold^d$ and for any $t \geq 1$ the time 
$\eta_t\left(\bx\right)$ is nothing 
but a $\mbox{Geometric}\left(\kappa\right)$
random variable. Thus,
\begin{eqnarray}
 &   & 
\Pbar^{\bzero}\left( \beta_t \leq t, \,\,\exists\,\,
                  \bx \in \Zbold^d
                  \mbox{\ s.t.\ }
                  \eta_{\beta_t}\left(\bx\right) \geq \beta_t
                  - \gamma_{\beta_t}\left(\bx\right) \right) 
                  \nonumber \\
 & \leq & 
C_0 \, \sum_{r=0}^{\infty} r^{d-1} \exp\left( - \lambda \, 
                              \left(L\left(t\right) \vee r \right) \right)
\nonumber \\
 &  = & 
C_0 \, \left( 
       \left( \sum_{r \leq L\left(t\right)} r^{d-1} \right) 
       \exp\left( - \lambda L\left(t\right) \right)
       +
       \sum_{r > L\left(t\right)} r^{d-1} \exp\left( - \lambda r \right)
       \right) \label{Equ:tau-Finite-Esti-2}
\end{eqnarray} 
where  $C_0 = C_0\left(d\right) > 0$  is such that 
an $L_1$ ball in $\Zbold^d$ of radius $r$ contains
less than 
$C_0 r^{d-1}$  points, 
and
$\lambda = - \log \left(1-\kappa\right)$. Indeed,
the first inequality in (\ref{Equ:tau-Finite-Esti-2})
follows from the observation that 
\[
\beta_t - \gamma_{\beta_t}\left(\bx\right) \geq 
L(t) \vee \vert \bx - X_{\beta_t - L(t)} \vert,
\]
with 
$\vert \,\cdot\, \vert$ denoting the $L_1$-norm on $\Zbold^d$, and
then computing the probability by conditioning on $X_{\beta_t - L(t)}$. 

Concerning the first term in (\ref{Equ:tau-Finite-Main-Esti}),  note that
\begin{eqnarray}
 &      & 
\Pbar^{\bzero} \left( \beta_t > t \right) \nonumber \\
 & \leq & 
\Pbar^{\bzero} \left( 
\mathop{\bigcap}\limits_{j=0}^{\lfloor \sfrac{t - L(t)}{L(t)} \rfloor}
\left[ \epsilon_{j L(t) + 1} =1, \epsilon_{j L(t) + 2} =1, \cdots,
       \epsilon_{ \left(j+1\right) L(t)} =1 \right]^c \right)
\nonumber \\
 & \leq &
\left( 1 -  \eps^{L(t)} \right)^
{\sfrac{t}{L(t)} - 2} \label{Equ:tau-Finite-Esti-1}
\end{eqnarray}

Using the choice of $L(t)$ and equations
(\ref{Equ:tau-Finite-Main-Esti}),
(\ref{Equ:tau-Finite-Esti-2}) and (\ref{Equ:tau-Finite-Esti-1}) 
one concludes that 
\begin{equation}
\Pbar^{\bzero} \left( \tau_1 > t \right) 
\leq
\frac{C_1}{t^{2 + \zeta}},
\label{Equ:tau-Finite-Final-Esti}
\end{equation}
where $C_1, \zeta > 0$ are some constants (depending on $\delta,d$). 
This completes the proof.
$\qed$

\begin{Remark}
\label{rem-190805} 
An inspection of the proof reveals that in fact,
\begin{equation}
\label{eq-190805a}
\forall \,\,\, \gamma'<\frac{\log(1-\kappa)}{\log \varepsilon}=:\gamma\,,
\quad \Ebar^{\bzero}\left[ \tau_1^{\gamma'} \right] < \infty\,.
\end{equation}
This will be useful when deriving the quenched invariance principle.
\end{Remark}


Define now a $\sigma$-algebra,
\begin{equation}
\HH_1 := \sigma \left( \tau_1 ; \left\{X_t\right\}_{0 \leq t \leq \tau_1} ;
\left\{ \left( \alpha_t\left(\cdot\right) \right)_{t \leq \tau_1} \right\}
\right).
\label{Equ:Defi-HH-1}
\end{equation}
The following is the most crucial lemma. 
\begin{Lemma}
\label{Lem:Renewal-Equ}
For any measurable sets $A, B, C$ we have,
\begin{align}
 & 
\Pbar^{\bzero} \left( 
\left\{ X_{\tau_1 + t} - X_{\tau_1} \right\}_{t \geq 0} \in A, \,
\left\{ \left( \omega_{\tau_1 + t}
\left( \cdot, \cdot \right) \right)_{t \geq 0} \right\} \in B, \,
\left\{ \left( \alpha_t\left(\cdot\right) \right)_{t \geq \tau_1 + 1} 
\right\} \in C 
\,\Bigg\vert\, \HH_1 \right) \nonumber \\
=&   
\Pbar^{\bzero} \left(
\left\{ X_t \right\}_{t \geq 0} \in A, \,
\left\{ \left( \omega_t\left( \cdot, \cdot \right) \right)_{t \geq 0}
\right\} \in B, \,
\left\{ \left( \alpha_t\left(\cdot\right) \right)_{t \geq 1} \right\} \in C 
\right). \label{Equ:Renewal-Equ}
\end{align}
\end{Lemma}
\proof Let $h$ be a $\HH_1$ measurable function. Write
$\bone_A := \bone \left( \left\{X_t - X_0 \right\}_{t \geq 0} \in A \right)$,
$\bone_B := \bone \left( \left\{ 
\left( \omega_t\left(\cdot, \cdot\right) \right)_{t \geq 0} \right\} \in B
\right)$ and
$\bone_C := \bone \left( 
\left\{ \left( \alpha_t \left( \cdot \right) \right)_{t \geq 1} \right\}
\in C \right)$. 
Note that for every $m \in \Nbold$ and $\bx \in \Zbold^d$, since
$h$ is $\HH_1$-measurable,
 there
exists a random variable $h_{\bx, m}$ which is measurable with respect to  
$\sigma \left( 
        \left\{ X_t \right\}_{t \geq 0}; 
        \left\{ \left( \alpha_t\left(\cdot\right) \right)_{t \leq m} \right\}
        \right)$,
such that $h  = h_{\bx, m}$
on the event $\left[ \tau_1 = m , \, X_{\tau_1} = \bx \right]$. 
Writing $\theta$ for the time shift,
\small
\begin{align*}
    & 
\Ebar^{\bzero}\left[ \bone_A \circ \theta^{\tau_1} \cdot
                  \bone_B \circ \theta^{\tau_1} \cdot
                  \bone_C \circ \theta^{\tau_1} \cdot h
              \right] \\
  = & 
\sum_{m=1}^{\infty} \, \sum_{\bx_m \in \Zbold^d} 
\Ebar^{\bzero}\left[ \bone_A \circ \theta^m \cdot
                  \bone_B \circ \theta^m \cdot
                  \bone_C \circ \theta^m \cdot 
                  \bone\left(\tau_1 = m \right) \cdot 
                  \bone\left(X_m = \bx_m \right) \cdot 
                  h_{\bx_m, m}
           \right] 
\end{align*}
\begin{align*}
  = &
\sum_{m=1}^{\infty} \, \sum_{\bx_m \in \Zbold^d} 
\Ebar^{\bzero}\left[ \bone_B \circ \theta^m \cdot 
                  \Ebb_{\omega}^{\bzero} \left[
                  \bone_A \circ \theta^m \cdot
                  \bone_C \circ \theta^m \cdot 
                  \bone\left(\tau_1 = m \right) \cdot 
                  \bone\left(X_m = \bx_m \right) \cdot 
                  h_{\bx_m, m}
           \right] \right] \\
  = &
\sum_{m=1}^{\infty}
\mathop{\sum_{\bx_t \in \Zbold^d}}\limits_{1 \leq t \leq m} \,
\mathop{\sum_{e_t \in W}}\limits_{1 \leq t \leq m}
\Ebar^{\bzero} \left[ 
\bone_B \circ \theta^m \cdot 
\Ebb_{\omega}^{\bzero} \left[
\bone_A \circ \theta^m \cdot
\bone_C \circ \theta^m \cdot 
\bone_{\left[\tau_1 = m \right]} \right.\right.\\
& \quad \quad \quad \quad \quad \quad \quad
 \quad \quad \quad \quad \quad \quad \quad
 \quad \quad \quad \quad
\left.\left.     \cdot
\bone_{\left[ X_t = \bx_t, 1 \leq t \leq m\right]} \cdot
\bone_{\left[ \epsilon_t= e_t, 1\leq t \leq m\right]} \cdot 
h_{\bx_m, m} 
\right]
\right] 
\end{align*}
\begin{align*}
  = &
\sum_{m=1}^{\infty}
\mathop{\sum_{\bx_t \in \Zbold^d}}\limits_{1 \leq t \leq m} \,
\mathop{\sum_{e_t \in W}}\limits_{1 \leq t \leq m}
\Ebar^{\bzero} \left[ 
\bone_B \circ \theta^m \cdot 
\Ebb_{\omega}^{\bzero} \left[
\Ebb_{\omega}^{\bx_m} \left[ \bone_A \cdot \bone_C \right] 
                      \circ \theta^m \cdot 
\bone_{\left[\tau_1 = m \right]}\right.\right. 
\\
& \quad \quad \quad \quad \quad \quad \quad
 \quad \quad \quad \quad \quad \quad \quad
 \quad \quad \quad \quad
\left.\left.
\cdot
\bone_{\left[ X_t = \bx_t, 1 \leq t \leq m\right]} \cdot
\bone_{\left[ \epsilon_t= e_t, 1\leq t \leq m\right]} \cdot 
h_{\bx_m, m} 
\right]
\right] \\
  = &
\sum_{m=1}^{\infty}
\mathop{\sum_{\bx_t \in \Zbold^d}}\limits_{1 \leq t \leq m} \,
\mathop{\sum_{e_t \in W}}\limits_{1 \leq t \leq m}
\Ebar^{\bzero} \left[ 
\Ebb_{\omega}^{\bzero} \left[
\Ebb_{\omega}^{\bx_m} \left[ \bone_A \cdot \bone_B \cdot \bone_C \right] 
                      \circ \theta^m \cdot 
\bone_{\left[\tau_1 = m \right]} \cdot
\bone_{\left[ X_t = \bx_t, 1 \leq t \leq m\right]} \right.\right.\\
& \quad \quad \quad \quad \quad \quad \quad
 \quad \quad \quad \quad \quad \quad \quad
 \quad \quad \quad \quad
 \quad \quad \quad \quad
  \quad \quad
\left.\left.
\cdot
\bone_{\left[ \epsilon_t= e_t, 1\leq t \leq m\right]} \cdot 
h_{\bx_m, m} 
\right]
\right] \\
  = &
\sum_{m=1}^{\infty}
\mathop{\sum_{\bx_t \in \Zbold^d}}\limits_{1 \leq t \leq m} \,
\mathop{\sum_{e_t \in W}}\limits_{1 \leq t \leq m}
\Ebar^{\bzero} \left[ 
\Ebb_{\omega}^{\bx_m} \left[ \bone_A \cdot \bone_B \cdot \bone_C \right] 
                      \circ \theta^m \cdot 
\bone_{\left[\tau_1 = m \right]} \cdot
\bone_{\left[ X_t = \bx_t, 1 \leq t \leq m\right]} \right.\\
& \quad \quad \quad \quad \quad \quad \quad
 \quad \quad \quad \quad \quad \quad \quad
 \quad \quad \quad \quad
 \quad \quad \quad \quad
  \quad \quad
\left.
\cdot
\bone_{\left[ \epsilon_t= e_t, 1\leq t \leq m\right]} \cdot 
h_{\bx_m, m} 
\right] \\
  = &
\Ebar^{\bzero}\left[\bone_A \cdot \bone_B \cdot \bone_C \right] \,
\sum_{m=1}^{\infty}
\mathop{\sum_{\bx_t \in \Zbold^d}}\limits_{1 \leq t \leq m} \,
\mathop{\sum_{e_t \in W}}\limits_{1 \leq t \leq m}
\Ebar^{\bzero} \left[
\bone_{\left[\tau_1 = m \right]} \cdot
\bone_{\left[ X_t = \bx_t, 1 \leq t \leq m\right]} \right.\\
& \quad \quad \quad \quad \quad \quad \quad
 \quad \quad \quad \quad \quad \quad \quad
 \quad \quad \quad \quad
 \quad \quad \quad \quad
  \quad \quad\left.
\cdot
\bone_{\left[ \epsilon_t = e_t, 1\leq t \leq m\right]} \cdot 
h_{\bx_m, m} 
\right],
\end{align*}
\normalsize
where in the fourth equality we use the Markov property of the random walk
given the environment and also the fact (\ref{Equ:alpha-independence}) 
with $s=m$; and the 
last equality uses the fact that 
on the event $\left[\tau_1 = m\right]$,
the environment chains at every site have gone
through ``regeneration'' before time $m$ and after the last ``proper''
visit of the walk to that site, and hence 
under the law $\Pbar^{\bzero}$, 
at time $m$ the environments at distinct sites are
independent,
$\pi$-distributed, and
independent of the $\epsilon$ and $\alpha$ coins  and the
walk till time $m$. 
Substituting in the above the whole sample space in place of $A, B$ and $C$,
we conclude that 
\begin{equation}
\Ebar^{\bzero}\left[ \bone_A \circ \theta^{\tau_1} \cdot
                  \bone_B \circ \theta^{\tau_1} \cdot
                  \bone_C \circ \theta^{\tau_1} \cdot h
           \right]
=
\Ebar^{\bzero}\left[\bone_A \cdot \bone_B \cdot \bone_C \right] \,
\Ebar^{\bzero}\left[ h \right],
\label{Equ:Final-Condi-Exp}
\end{equation}
concluding the proof of the lemma. $\qed$

Consider now $\tau_1$ as a function of of
$\left( \left(X_t\right)_{t \geq 0}; 
\left\{ \left( \omega_t\left(\cdot, \cdot \right) \right)_{t \geq 0} \right\}; 
\left\{ \left(\alpha_t\left(\cdot\right) \right)_{t \geq 1} \right\} \right)$, 
and set 
\begin{equation}
\tau_{n+1} := \tau_n + 
\tau_1 \left( \left(X_{\tau_n + t}\right)_{t \geq 0}; 
\left\{ \left( \omega_{\tau_n+ t}
\left(\cdot, \cdot \right) \right)_{t \geq 0} \right\}; 
\left\{ \left(\alpha_{\tau_n+t}\left(\cdot\right) \right)_{t \geq 1} 
\right\} \right),
\label{Equ:Defi-tau-n}
\end{equation}
with $\tau_{n+1} = \infty$ on the event $\left[ \tau_n = \infty \right]$. 
The following lemma gives the renewal sequence described earlier.
\begin{Lemma}
\label{Lem:Renewal}
$\Pbar^{\bzero} \left(\tau_n < \infty \right) = 1 $, 
for all $n \geq 1$. Moreover 
the sequence of random vectors  
$\left\{ \left( \tau_{n+1} - \tau_n , X_{\tau_{n+1}} - X_{\tau_n} \right)
\right\}_{n \geq 0}$, where $\tau_0 = 0$, are i.i.d. under the law
$\Pbar^{\bzero}$. 
\end{Lemma} 
\proof Define
\begin{equation}
\HH_n := \sigma \left( \tau_1, \tau_2, \ldots, \tau_n ;
\left(X_t\right)_{0 \leq t \leq \tau_n} ;
\left\{ \left( \alpha_t\left(\bx\right) \right)_{t \leq \tau_n} \right\}
\right),
\label{Equ:Defi-HH-n}
\end{equation}
then an obvious rerun of the proof of Lemma \ref{Lem:Renewal} yields that
for measurable sets $A, B$ and $C$, 
\begin{align}
  \ &
\Pbar^{\bzero} \left( 
\left\{ X_{\tau_n + t} - X_{\tau_n} \right\}_{t \geq 0} \in A, \,
\left\{ \left( \omega_{\tau_n + t}
\left( \cdot, \cdot \right) \right)_{t \geq 0} \right\} \in B, \,
\left\{ \left( \alpha_{\tau_n + t}\left(\cdot\right) \right)_{t \geq 1} 
\right\} \in C 
\,\Bigg\vert\, \HH_n \right) \nonumber \\
 = & 
\Pbar^{\bzero} \left(
\left\{ X_t \right\}_{t \geq 0} \in A, \,
\left\{ \left( \omega_t\left( \cdot, \cdot \right) \right)_{t \geq 0}
\right\} \in B \,
\left\{ \left( \alpha_t\left(\cdot\right) \right)_{t \geq 1} \right\} \in C 
\right). \label{Equ:Re-Renewal-Equ}
\end{align}
So first of all we get from Proposition \ref{Prop:tau-Finite} that 
$\Pbar^{\bzero} \left( \tau_n < \infty \right) = 1$, for all $n \geq 1$, 
and also under $\Pbar^{\bzero}$
\[
\left( \tau_1, X_{\tau_1} \right), 
\left( \tau_2 - \tau_1, X_{\tau_1} - X_{\tau_1} \right), 
\ldots, 
\left( \tau_{n+1} - \tau_n, X_{\tau_{n+1}} - X_{\tau_n} \right),
\ldots
\]
are i.i.d. of random vectors. $\qed$

\section{Proofs of the main results}
\label{Sec:Proofs}

\subsection{Proof of Theorem \ref{Thm:SLLN}}
\label{Subsec:Proof-of-Them-SLLN}

Fix $\ba \in \Rbold^d$. 
From Proposition \ref{Prop:tau-Finite} we get that 
$\Ebar^{\bzero}\left[ \tau_1 \right] < \infty$. Since the random walk
$\left\{ X_t \right\}_{t \geq 0}$ has bounded increments,  it follows that 
$\Ebar^{\bzero}\left[ \ba \cdot X_{\tau_1} \right] < \infty$ as well. 
Let $Y_n := \ba \cdot \left( X_{\tau_n} - X_{\tau_{n-1}} \right)$, for 
$n \geq 1$, taking $\tau_0 = 0$. Using Lemma \ref{Lem:Renewal}
we conclude
\begin{equation}
\frac{\tau_n}{n} \mathop{\longrightarrow}\limits_{n\to\infty} 
\Ebar^{\bzero}\left[ \tau_1 \right] \,\,\,\, \mbox{a.s.\ \ }
\left[\Pbar^{\bzero}\right], 
\label{Equ:Limit-tau-n}
\end{equation}
and
\begin{equation}
\lim_{n \rightarrow \infty} \frac{1}{n} \sum_{k=1}^n Y_k
\mathop{\longrightarrow }\limits_{n\to\infty}
\Ebar^{\bzero} \left[ \ba \cdot X_{\tau_1} \right] \,\,\,\, 
\mbox{a.s.\ \ } \left[\Pbar^{\bzero}\right]. 
\label{Equ:Limit-Sum-Y}
\end{equation}

We can find a (possibly random) sequence of numbers 
$\left\{ k_n \right\}_{n \geq 1}$ increasing to 
$\infty$ such that for all $n \geq 1$ we have 
$\tau_{k_n} \leq n < \tau_{k_n + 1}$. Then
from (\ref{Equ:Limit-tau-n})
we get
\begin{equation}
\frac{n}{k_n} \mathop{\longrightarrow}\limits_{n\to\infty}
\Ebar^{\bzero} \left[ \tau_1 \right] \,\,\,\, 
\mbox{a.s.\ \ } \left[\Pbar^{\bzero}\right]. 
\label{Equ:Limit-k-n}
\end{equation}
Since the increments of the random walk are bounded, we have
for any $n \geq 1$  that
\begin{eqnarray}
\Big\vert \ba \cdot\left(X_n - X_{\tau_{k_n}}\right) \Big\vert
 & \leq & \Vert \ba \Vert_2 \, \left( n - \tau_{k_n} \right) 
   \nonumber \\
 & \leq & \Vert \ba \Vert_2 \, 
          \left( \tau_{k_n + 1} - \tau_{k_n} \right)\,.
   \label{Equ:SLLN-Err-Esti}
\end{eqnarray}
On the other hand, since 
by  Lemma \ref{Lem:Renewal} the random variables
$\left( \tau_n - \tau_{n-1} \right)_{n \geq 1}$ are identically
distributed and of finite mean,
one gets for any $\delta>0$ that
$$ 
\sum_{n=0}^\infty\Pbar^{\bzero}\left(\tau_{n+1}-\tau_n>\delta n\right)
= 
\sum_{n=0}^\infty\Pbar^{\bzero}\left(\tau_{1}>\delta n\right)
<\infty\,.$$
It follows from an application of the Borel-Cantelli Lemma that
\begin{equation}
\frac{\tau_{n+1} - \tau_n}{n} \mathop{\longrightarrow}\limits_{n\to\infty} 0
\,\,\,\, \mbox{a.s.\ \ } \left[\Pbar^{\bzero}\right]. 
\label{Equ:Apply-BC-tau}
\end{equation}
So using (\ref{Equ:SLLN-Err-Esti}) and (\ref{Equ:Limit-k-n}) we get 
\begin{equation}
\frac{\Big\vert \ba \cdot\left(X_n - X_{\tau_{k_n}}\right) \Big\vert}{n}
\mathop{\longrightarrow}\limits_{n\to\infty} 0
\,\,\,\, \mbox{a.s.\ \ } \left[\Pbar^{\bzero}\right]
\label{Equ:SLLN-Err-Limit}
\end{equation}
This together with (\ref{Equ:Limit-Sum-Y}) and (\ref{Equ:Limit-k-n}) 
gives 
\begin{equation}
\frac{\ba \cdot X_n}{n} \mathop{\longrightarrow}\limits_{n\to\infty}
\frac{\Ebar^{\bzero}\left[ \ba \cdot X_{\tau_1} \right]}
     {\Ebar^{\bzero}\left[\tau_1\right]}
\,\,\,\, \mbox{a.s.\ \ } \left[\Pbar^{\bzero}\right]
\label{Equ:Final-Limit-a}
\end{equation}
Finally taking $\ba = \be_i$ for $i = 1, 2, \ldots, d$ we conclude
\begin{equation}
\frac{X_n}{n} \mathop{\longrightarrow}\limits_{n\to\infty}
\frac{\Ebar^{\bzero}\left[ X_{\tau_1} \right]}
     {\Ebar^{\bzero}\left[\tau_1\right]}
\,\,\,\, \mbox{a.s.\ \ } \left[\Pbar^{\bzero}\right]
\label{Equ:Final-Limit}
\end{equation}
which completes the proof. $\qed$
\subsection{Proof of Theorem \ref{Thm:CLT}}
\label{Subsec:Proof-of-Thm-CLT}
Fix $\ba \in \Rbold^d$, let $\left(Y_n\right)_{n \geq 1}$ be as defined
above. Put 
$\Ybar_n := Y_n - \left(\tau_n - \tau_{n-1}\right) \ba \cdot \bv$, where
$\bv = \sfrac{\Ebar^{\bzero}\left[ X_{\tau_1} \right]}
             {\Ebar^{\bzero}\left[\tau_1\right]}$.
Let $S_n := \Ybar_1 + \Ybar_2 + \cdots + \Ybar_n$, for $n \geq 1$. 
By Proposition \ref{Prop:tau-Finite}, $\tau_1$ has finite second moment and,
due to the boundedness of the increments
of the random walk $\left( X_t \right)_{t \geq 0}$,
so does $X_{\tau_1}$. Further, by definition, 
 $\Ebar^{\bzero}\left[ S_n \right]=0$. 
Thus by Lemma \ref{Lem:Renewal} and Donsker's invariance principle,
see e.g. \cite[Theorem 14.1]{Bill99},
we have
\begin{equation}
\frac{S_{\lfloor nt \rfloor
}}{\sigma_a \sqrt{n}} \mathop{\longrightarrow}\limits^d_{n\to\infty} 
\mbox{BM}\left(1\right),
\label{Equ:Standard-CLT}
\end{equation}
where $\mbox{BM}(1):=\mbox{BM}_1(1)$ and
$\sigma_a^2 := \Ebar^{\bzero}\left[ \Ybar_1^2 \right]>0$, where
the last inequality is due to
Condition (A2) and (\ref{Equ:q-non-trivial}).

Put next  
$m_n := n / \Ebar^{\bzero}\left[\tau_1\right] $.
Since $n\mapsto  m_n $ is a deterministic scaling,
it follows from
(\ref{Equ:Standard-CLT}) that
\begin{equation}
\frac{S_{\lfloor m_n t \rfloor}}{\sigma_a \sqrt{m_n}} 
\mathop{\longrightarrow}\limits^d_{n\to\infty} 
\mbox{BM}\left(1\right).
\label{Equ:CLT-in-m-n}
\end{equation}

Let $\left(k_n\right)_{n \geq 1}$ be as before. 
Our next step is to prove the analogue of 
(\ref{Equ:Standard-CLT}) with $m_nt$ replaced by $k_{\lfloor nt \rfloor}$.
A consequence of 
(\ref{Equ:Limit-k-n}) is that for any fixed $T<\infty$,
\begin{equation}
\sup_{t\leq T}
\left|\frac{ k_{\lfloor nt\rfloor} }{n}
-\frac{\lfloor m_n t\rfloor }{n}\right| 
\mathop{\longrightarrow}\limits_{n\to\infty} 0\,,
\,\,\,\, \mbox{a.s.\ \ } \left[\Pbar^{\bzero}\right]
\label{Equ:Limit-m-n}
\end{equation}
Thus,
\begin{equation}
\frac{S_{ k_{\lfloor n t\rfloor}}}{\sigma_a \sqrt{m_n}}
\mathop{\longrightarrow}\limits^d
\mbox{BM}\left(1\right).
\label{Equ:CLT-in-k-n}
\end{equation}
Note that for an appropriate $C=C(\ba,d)$, using (\ref{Equ:SLLN-Err-Esti}),
\begin{equation}
\label{eq-180805a}
\sup_{t\leq T}\left|
\frac{ S_{k_{\lfloor nt\rfloor}}- \ba \cdot X_{\lfloor nt\rfloor}-nt \, \ba \cdot 
\bv}{\sqrt{n}}
\right|\leq
C \max_{0\leq i\leq k_{\lfloor nT\rfloor}} \frac{\tau_{i+1}-\tau_i}{\sqrt{n}}\,.
\end{equation}
Since 
$k_{\lfloor nT\rfloor}\leq \lfloor nT \rfloor$,
we have that for any $\delta>0$,
\begin{equation}
\label{eq-180805b}
\Pbar^{\bzero}
\left( \max_{0\leq i\leq k_{\lfloor nT\rfloor}} 
\frac{\tau_{i+1}-\tau_i}{\sqrt{n}}
 >\delta \right)
\leq 
\sum_{i=1}^{\lfloor nT \rfloor} 
\Pbar^{\bzero}(\tau_1>\delta \sqrt n)\,.
\end{equation}
Note that, since
$
\Ebar^{\bzero}\left[\tau_1^2\right]<\infty$, one has that
$$ \sum_{i=1}^{\infty} 
\Pbar^{\bzero}
(\tau_1>\frac{\delta \sqrt{i}}{\sqrt T})=
\sum_{i=1}^\infty 
\Pbar^{\bzero}
(\tau_1^2>\frac{\delta^2 i}{T})<\infty\,.$$
Hence, for each $\delta_1>0$ there is a deterministic constant
 $I_{\delta_1}$ depending 
on $d,\delta,T,\delta_1$   such that
$$ \sum_{i=I_{\delta_1}}^{\infty} 
\Pbar^{\bzero}
(\tau_1>\frac{\delta \sqrt{i}}{\sqrt{T}})<
\delta_1\,.$$
Therefore,
\begin{eqnarray*}
&  & \limsup_{n\to\infty}\sum_{i=1}^{\lfloor nT \rfloor} 
\Pbar^{\bzero}(\tau_1>\delta \sqrt n) \\
& \leq & \limsup_{n\to\infty} \left( \sum_{i=1}^{I_{\delta_1}} 
\Pbar^{\bzero}
(\tau_1>\delta\sqrt{n})+
\sum_{i=I_{\delta_1}+1}^\infty 
\Pbar^{\bzero}
(\tau_1>\frac{\delta \sqrt{i}}{ \sqrt{T}}) \right)
\leq \delta_1\,.
\end{eqnarray*}
$\delta_1$ being arbitrary, one concludes from the last limit and 
(\ref{eq-180805b}) that
$$
\Pbar^{\bzero}
\left( \max_{0\leq i\leq k_{\lfloor nT\rfloor}} 
\frac{\tau_{i+1}-\tau_i}{\sqrt{n}}
 >\delta \right)
\mathop{\longrightarrow}\limits_{n\to\infty} 0\,.$$
Together with
(\ref{Equ:CLT-in-k-n}) and 
(\ref{eq-180805a}), this implies that
for every $\ba \in \Rbold^d\setminus 0$ we have
\begin{equation}
\frac{\ba\cdot\left(X_{ {\lfloor n t\rfloor}}-nt \, \bv\right)}
{\sigma_a \sqrt{m_n}} \mathop{\longrightarrow}\limits^d \mbox{BM}\left(1\right).
\label{eq-180805c}
\end{equation}
Since $\ba$ is arbitrary, this completes the proof of the theorem,
with
\begin{equation}
\label{eq-250805aa}
\Sigma := \frac{\Var_{\Pbar^{\bzero}}\left( X_{\tau_1} - \tau_1 \, \bv \right)}
{\Ebar^{\bzero}\left[\tau_1\right]}.\end{equation}
 $\qed$

\subsection{Proof of Theorem \ref{Thm:CLTquenched}}
\label{Subsec:Proof-of-Thm-CLTquenched}
Our argument is based on the technique introduced by 
Bolthausen and Sznitman in \cite{BS}, as developed in
\cite{BSZ}. Let $B^n_t=(X_{\lfloor nt \rfloor}-nt \bv)/\sqrt{n}$,
and let $\BBB^n_t$ denote the polygonal interpolation 
of $(k/n)\to B^n_{k/n}$. 
Consider the space 
$\CC_T:=C([0,T],\Rbold^d)$ of continuous $\Rbold^d$-valued functions 
on $[0,T]$, endowed with the distance
$d_T(u,u')=\sup_{t\leq T}|u(t)-u'(t)|\wedge 1$,
By \cite[Lemma 4.1]{BS},
Theorem \ref{Thm:CLTquenched} follows from Theorem \ref{Thm:CLT}
once we show that for all bounded Lipschitz function $F$ on
$\CC_T$ and $b\in (1,2]$,
\begin{equation}
\label{eq-200805b}
\sum_{m=1}^{\infty} \Var_{\bP^{\pi}}\left(\bE^0_\omega
[F(\BBB^{\lfloor b^m
\rfloor})]\right)<\infty\,.
\end{equation}
In order to prove (\ref{eq-200805b}), we now follow the
approach of \cite{BSZ}. We construct the environment
$\omega$ using the variables $(\alpha_t(\bx))_{t\geq 1, 
\bx\in \Zbold^d}$ as described in Section 
\ref{sec-regeneration}. We next construct two independent 
sequences of i.i.d. Bernoulli($\varepsilon$) random variables,
that we denote by $(\varepsilon_t^{(1)})_{t\geq 2}$ and
$(\varepsilon_t^{(2)})_{t\geq 1}$. Given these
sequences and the environment $\omega$, we construct two independent
copies of the random walk, denoted $(X_t^{(1)})_{t\geq 1}$ and
$(X_t^{(2)})_{t\geq 1}$, following the recipe of Section 
\ref{sec-regeneration} (we of course use the sequence
$\varepsilon^{(j)}$ to construct $X^{(j)}$, for $j=1,2$),
and introduce the respective linear interpolations
$(\BBB^{n,(j)}_t)_{t\geq 0}$.
It is then clear that (\ref{eq-200805b}) is equivalent to 
\begin{align}
\label{eq-200805c}
& \sum_m \left(
\bP^{\pi}\otimes \bQ_\varepsilon
\otimes
\bQ_\varepsilon
\left(\Pbar_{\omega, {\mathbf \epsilon}^{(1)}}^{\bzero} 
\otimes \Pbar_{\omega, {\mathbf \epsilon}^{(2)}}^{\bzero} 
\left[F(\BBB^{\lfloor b^m
\rfloor, (1)})
F(\BBB^{\lfloor b^m
\rfloor, (2)})\right]\right)\right. \\
&\left.\quad \quad -
\bP^{\pi}
\otimes \bP^{\pi}\otimes 
\bQ_\varepsilon\otimes
\bQ_\varepsilon
\left(\Pbar_{\omega^{(1)}, {\mathbf \epsilon}^{(1)}}^{\bzero} 
\otimes \Pbar_{\omega^{(2)}, {\mathbf \epsilon}^{(2)}}^{\bzero} 
\left[F(\BBB^{\lfloor b^m
\rfloor, (1)})
F(\BBB^{\lfloor b^m
\rfloor, (2)})\right]\right)
\right)\nonumber \\
& 
\quad\quad \quad 
\quad\quad \quad 
\quad\quad \quad 
\quad\quad \quad 
\quad\quad \quad 
\quad\quad \quad 
\quad\quad \quad 
\quad\quad \quad 
\quad\quad \quad 
\quad\quad \quad <\infty\, ,
\nonumber \end{align}
where for any probability measure $\bP$ and a measurable function $f$
by $\bP\left(f\right)$ we mean $\bE_{\bP}\left[f\right]$. 
In the sequel, we write
$\Pbar_{\omega, \omega', {\mathbf \epsilon},{\mathbf \epsilon}'}
^{\bzero}$ for  
$\Pbar_{\omega, {\mathbf \epsilon}}^{\bzero} 
\otimes \Pbar_{\omega', {\mathbf \epsilon}'}^{\bzero} $.
Recalling the constant $\gamma$ from (\ref{eq-200805na}),
we next choose
 constants $\theta,\theta',\theta'',\mu,\alpha$ satisfying the following
conditions:
\begin{align}
& 0< \theta<1, \,\, 2/\gamma<\theta'<\theta/2, 
 \,\, \theta>2(\theta'+1/(d-1)), \,\, \theta''<\theta',\,
(\theta'-\theta'')\gamma>1, 
\nonumber \\
& 0<\mu<1/2, \,\, 1/2>\alpha >  (1/\theta'+1)/\gamma, \nonumber \\
& \theta''((d-1)-4-2\alpha (d-1)) 
= \theta''(d - 5 - 2\alpha d + 2 \alpha) > 1\,.
\label{cond-constants}
\end{align}
It is not hard to verify that
the assumptions of 
Theorem \ref{Thm:CLTquenched}
 imply that such constants can be found
(indeed, verify, using that $d>7$ 
and $d-1>4/(1-6/\gamma)$ which implies
that $\gamma>6(d-1)/(d-5)$,
 that taking $\theta'=(d-3)/2(d-1)$,
$\alpha=1/\gamma+1/\gamma\theta'$  and
$\theta''= \theta' (\gamma(d-3)-2(d-1))/\gamma (d-3)$ satisfies 
the constraints except for the last one with equality,
 and the last one with
inequality due to (\ref{eq-200805a}), and argue by 
continuity).
Fix then an integer $m$.
We let $\bar \beta_m$ denote the first time
there is a run of length $\lfloor m^{\theta} \rfloor$ of non-zero
$\epsilon$-coin tosses of both types ending at it, that is
$$ 
\bar \beta_m=\inf\{s\geq \lfloor m^{\theta} \rfloor\,\Big\vert\;
\epsilon_s^{(1)}= \epsilon_s^{(2)}=
\epsilon_{s-1}^{(1)}= \epsilon_{s-1}^{(2)}=
\ldots =\epsilon_{s-\lfloor m^{\theta} \rfloor+1}^{(1)}= \epsilon_{s-\lfloor m^{\theta} \rfloor+1}^{(2)}=
1\}\,.$$
For any $\bx\in \Zbold^d$, set $D^{(j)}_m(\bx)
=|\bx - X^{(j)}_{\bar \beta_m}|$.
Setting $n=\lfloor b ^m \rfloor$,
define next the events
$$ G_m=\{\bar\beta_m\leq n^\mu,  
X^{(1)}_{[\bar\beta_m,\infty)}\cap X^{(2)}_{[\bar \beta_m,\infty)}
=\emptyset\}\,,$$
and for $j=1,2$,
$$R_m^{(j)}=\{\forall \bx\in \Zbold^d, \exists t\in  
[\bar\beta_m-\lfloor m^{\theta} \rfloor, \bar\beta_m+D^{(j)}_m(\bx)] 
\,\mbox{\rm such that}\,
\, \alpha_t(\bx)=1\}\,.$$
Finally, set $\hat G_m=G_m\cap R_m^{(1)}\cap R_m^{(2)}$.
The crucial element of the proof of Theorem \ref{Thm:CLTquenched}
is contained in the following lemma, whose proof is postponed.
\begin{Lemma}
\label{lem-quenched}
Under the assumptions of Theorem \ref{Thm:CLTquenched},
the following estimates hold.
\begin{align}
\label{eq-200805n1}
&\sum_{m=1}^\infty
\bP^{\pi}\otimes \bQ_\varepsilon
\otimes
\bQ_\varepsilon
\left(\Pbar_{\omega, \omega,{\mathbf \epsilon}^{(1)},
{\mathbf \epsilon}^{(2)}}^{\bzero} 
[(\hat G_m)^c]\right)<\infty\,.\\
\label{eq-200805n2}
&\sum_{m=1}^\infty
\bP^{\pi}\otimes 
\bP^{\pi}\otimes 
\bQ_\varepsilon
\otimes
\bQ_\varepsilon
\left(\Pbar_{\omega^{(1)}, \omega^{(2)},
{\mathbf \epsilon}^{(1)} 
, {\mathbf \epsilon}^{(2)}}^{\bzero} 
[(\hat G_m)^c]\right)<\infty\,.
\end{align}
\end{Lemma}
Equipped with Lemma \ref{lem-quenched}, let us complete
the proof (\ref{eq-200805c}). Indeed, for any integer $m$,
\begin{align}
\label{eq-200805d}
\Delta_m:= \nonumber \\
 & \Pbar_{\omega^{(1)}, 
\omega^{(1)}, 
{\mathbf \epsilon}^{(1)}, 
{\mathbf \epsilon}^{(2)}}^{\bzero} 
\left[F(\BBB^{n
, (1)})
F(\BBB^{n
, (2)})\right]
-\Pbar_{\omega^{(1)}, 
\omega^{(2)}, 
{\mathbf \epsilon}^{(1)}
{\mathbf \epsilon}^{(2)}}^{\bzero} 
\left[F(\BBB^{n
, (1)})
F(\BBB^{n
, (2)})\right]\nonumber \\
=&
\Pbar_{\omega^{(1)}, 
\omega^{(1)}, 
{\mathbf \epsilon}^{(1)}, 
{\mathbf \epsilon}^{(2)}}^{\bzero} 
\left[F(\BBB^{n
, (1)})
F(\BBB^{n
, (2)}); \hat G_m\right]\nonumber \\
&\quad 
-\Pbar_{\omega^{(1)}, 
\omega^{(2)}, 
{\mathbf \epsilon}^{(1)}
{\mathbf \epsilon}^{(2)}}^{\bzero} 
\left[F(\BBB^{n
, (1)})
F(\BBB^{n
, (2)});\hat G_m\right]+ d_m 
\end{align}
where
$$ \sum_{m=1}^\infty \bP^{\pi}\otimes 
\bP^{\pi}\otimes 
\bQ_\varepsilon
\otimes
\bQ_\varepsilon( \vert d_m \vert )<\infty\,.$$
Let $\hat{\BBB}^{n,(j)}_t=
\BBB^{n,(j)}_{t+\sfrac{\bar\beta_m}{n}}
-\BBB^{n,(j)}_{\sfrac{\bar\beta_m}{n}}$ for $j=1,2$.
Recall that on $\hat G_m$ one has $\bar \beta_m<n^\mu$
and hence, using the Lipschitz property of $F$, 
it holds a.s. with respect to $\Pbar^{\bzero}$ on $\hat G_m$ that
$$
|F(\BBB^{n
, (j)})-F(\hat{\BBB}^{n,(j)})|\leq 2n^{\mu-1/2}\,.$$
Substituting in (\ref{eq-200805d}), we get
\begin{align}
\label{eq-200805e}
\Delta_m=&
\Pbar_{\omega^{(1)}, 
\omega^{(1)}, 
{\mathbf \epsilon}^{(1)}, 
{\mathbf \epsilon}^{(2)}}^{\bzero} 
\left[F(\hat{\BBB}^{n
, (1)})
F(\hat{\BBB}^{n
, (2)}); \hat G_m\right]\nonumber \\
&\quad 
-\Pbar_{\omega^{(1)}, 
\omega^{(2)}, 
{\mathbf \epsilon}^{(1)}
{\mathbf \epsilon}^{(2)}}^{\bzero} 
\left[F(\hat{\BBB}^{n
, (1)})
F(\hat{\BBB}^{n
, (2)});\hat G_m\right]+ e_m
\end{align} 
where 
\begin{equation}
\label{eq-200805f}
\sum_{m=1}^\infty 
\bP^{\pi}\otimes 
\bP^{\pi}\otimes 
\bQ_\varepsilon
\otimes
\bQ_\varepsilon( \vert e_m \vert )<\infty\,.
\end{equation}
Conditioning on $X_{\bar \beta_m}^{(j)}$, $j=1,2$, one observes that
\begin{align}
&\bP^{\pi}\otimes 
\bP^{\pi}\otimes 
\bQ_\varepsilon
\otimes
\bQ_\varepsilon\left(
\Pbar_{\omega^{(1)}, 
\omega^{(1)}, 
{\mathbf \epsilon}^{(1)}, 
{\mathbf \epsilon}^{(2)}}^{\bzero} 
\left[F(\hat{\BBB}^{n
, (1)})
F(\hat{\BBB}^{n
, (2)}); \hat G_m\right]\right.\nonumber \\
&\left.\quad 
-\Pbar_{\omega^{(1)}, 
\omega^{(2)}, 
{\mathbf \epsilon}^{(1)}
{\mathbf \epsilon}^{(2)}}^{\bzero} 
\left[F(\hat{\BBB}^{n
, (1)})
F(\hat{\BBB}^{n
, (2)});\hat G_m\right]\right)=0\,.
\end{align} 
Together with (\ref{eq-200805e}) and (\ref{eq-200805f}), one concludes 
that
$$\sum_{m=1}^\infty 
\bP^{\pi}\otimes 
\bP^{\pi}\otimes 
\bQ_\varepsilon
\otimes
\bQ_\varepsilon(\Delta_m)<\infty\,,
$$
as claimed.
$\qed$

\noindent
{\bf Proof of lemma \ref{lem-quenched}}
We begin by considering the event $\{\bar\beta_m>n^\mu\}$. By
the independence of the i.i.d. sequences $\epsilon^{(j)}$, we get
the estimate
\begin{eqnarray}
\label{eq-lem200805a}
\bQ_\varepsilon
\otimes
\bQ_\varepsilon(\bar\beta_m>n^\mu)
&\leq &(1-\varepsilon^{2\lfloor 
m^{\theta} \rfloor})^{n^\mu/\lfloor m^{\theta} \rfloor - 2}\nonumber \\
&\leq &
e^{-n^\mu \varepsilon^{2\lfloor m^{\theta} 
\rfloor}/\lfloor m^{\theta} \rfloor
+ 2 \eps^{2 \lfloor m^{\theta} \rfloor}}
\leq c_1 e^{-e^{c_2 m}}
\,,
\end{eqnarray}
for appropriate constants $c_1$, $c_2$, where the last 
estimate used that $\theta<1$. 

We next consider the event $(R_m^{(j)})^c\cap\{\bar\beta_m\leq 
n^\mu\}$. Decomposing according to the distance from $X_{\bar\beta_m}
^{(j)}$, as in the proof of Proposition  
\ref{Prop:tau-Finite} (see (\ref{Equ:tau-Finite-Esti-2})), one
 gets for some deterministic constants $c_3$,
$c_4$ (that may depend on the choice of parameters)
\begin{equation}
\label{eq-lem200805b}
\bP^{\pi}\otimes 
\bQ_\varepsilon
((R_m^{(j)})^c\cap\{\bar\beta_m\leq 
n^\mu\})\leq 
c_3\sum_{r=1}^\infty r^{d-1}(1-\kappa)^{r+m^\theta}
\leq  c_4 (1-\kappa)^{m^\theta}\,.
\end{equation}

We next turn to the crucial estimate of the probability of
non-intersection after $\bar \beta_m$. For this we consider the
cases $\bv = \bzero$ and $\bv \neq \bzero$ separately. 

First suppose $\bv \neq \bzero$. 
We start by showing that at time $\bar\beta_m$, the two walkers are not
likely to be too close in a $(d-1)$-dimensional sub-space. 
More precisely, 
let $V^{\bot}$ be the $(d-1)$ dimensional sub-space of $\Rbold^d$
which is orthogonal to the vector $\bv$. Let $\Lbold^{d-1}$ be the 
$(d-1)$-dimensional 
lattice which is the projection of $\Zbold^d$ into $V^{\bot}$. 
Let $P^{\bx}_q$ denote the law 
of a homogeneous Markov chain on $\Zbold^d$ with transition 
probabilities determined by $q$ starting at $\bx \in \Zbold^d$
and let $(M^{\tilde\bx}_t)_{t\geq 0}$
denote the projection of the
associated walk into the lattice $\Lbold^{d-1}$, where
$\tilde{\bx}$ is the projection of the point $\bx$ into the lattice
$\Lbold^{d-1}$. 
In particular,
$P^{\bx}_q(M^{\tilde{\bx}}_0=\tilde{\bx})=1$. We use 
$M^{\tilde{\bx}},\bar M^{\tilde{\by}}$
to denote independent copies of such walks starting at $\bx,\by$,
respectively. For $j=1,2$, let
$\widetilde{X}_{\bar\beta_m}^{(j)}$ denote the 
projection of the vector $X_{\bar\beta_m}^{(j)}$
into $V^{\bot}$.
Then, with
$A_m=\{|\widetilde{X}_{\bar\beta_m}^{(1)}-
\widetilde{X}_{\bar \beta_m}^{(2)}|\geq \lfloor m^{\theta'} \rfloor\}$,
we have, 
\begin{align}
\label{eq-210805a}
&\bP^{\pi}\otimes 
\bP^{\pi}\otimes 
\bQ_\varepsilon
\otimes
\bQ_\varepsilon\left(
\Pbar_{\omega^{(1)}, 
\omega^{(1)}, 
{\mathbf \epsilon}^{(1)}, 
{\mathbf \epsilon}^{(2)}}^{\bzero} 
(A_m^c)\right)\nonumber \\
\leq & \sup_{\bx,\by\in \Zbold^d}
P^{\bx}_q\otimes P^{\by}_q
\left(|M^{\tilde{\bx}}_{\lfloor m^{\theta} \rfloor}
- M^{\tilde{\by}}_{\lfloor m^{\theta} \rfloor}| < m^{\theta'}\right) \nonumber \\
= & \sup_{\bx, \by \in \Zbold^d} \,
P_q^{\by} \left( P_q^{\bx} \left( \vert M_{\lfloor m^{\theta} \rfloor}^{\tilde{\bx}} 
- M_{\lfloor m^{\theta} \rfloor}^{\tilde{\by}} \vert < m^{\theta'} \,\Big\vert\, 
M_{\lfloor m^{\theta} \rfloor}^{\tilde{\by}} \,\right) \right) \nonumber \\
\leq & \, C m^{\theta' (d-1)} \,\,
\sup_{\tilde{\bx}, \tilde{\bz} \in \Lbold^{d-1}} P_q^{\bx} \left( 
M_{\lfloor m^{\theta} \rfloor}^{\tilde{\bx}} = \tilde{\bz} \right) \nonumber \\
\leq & \,\, C m^{\theta' (d-1)} \frac{1}{m^{\theta (d-1)/2}}=
Cm^{-\left(\frac{\theta}{2}-\theta'\right)(d-1)} 
\end{align}
where $C$ is some deterministic constant and the last inequality
is due to the local limit theorem for
lattice distribution under $P^{\bx}_q$, see \cite[Theorem 22.1]{BhaRao76}. 

Letting
$$B_m =A_m\cap \{\bar\beta_m \leq n^\mu\}
\cap R_m^{(1)}\cap R_m^{(2)},$$ and
repeating the estimate (\ref{eq-210805a})
without
change for 
$\Pbar_{\omega^{(1)}, 
\omega^{(2)}, 
{\mathbf \epsilon}^{(1)}, 
{\mathbf \epsilon}^{(2)}}^{\bzero}$ replacing 
$\Pbar_{\omega^{(1)}, 
\omega^{(1)}, 
{\mathbf \epsilon}^{(1)}, 
{\mathbf \epsilon}^{(2)}}^{\bzero} $,
and using that
$(d-1)(\theta/2-\theta')>1$, we conclude that
\begin{align}
\label{eq-210805b}
&\max\left(
\sum_{m=1}^\infty \bP^{\pi}\otimes 
\bP^{\pi}\otimes 
\bQ_\varepsilon
\otimes
\bQ_\varepsilon\left(
\Pbar_{\omega^{(1)}, 
\omega^{(1)}, 
{\mathbf \epsilon}^{(1)}, 
{\mathbf \epsilon}^{(2)}}^{\bzero} 
(B_m^c)\right),\right.\\
&
\quad\quad \quad \quad \left.\sum_{m=1}^\infty \bP^{\pi}\otimes 
\bP^{\pi}\otimes 
\bQ_\varepsilon
\otimes
\bQ_\varepsilon\left(
\Pbar_{\omega^{(1)}, 
\omega^{(2)}, 
{\mathbf \epsilon}^{(1)}, 
{\mathbf \epsilon}^{(2)}}^{\bzero} 
(B_m^c)\right)\right)<\infty
\nonumber \end{align}
Next, for $j=1,2$,
we construct, using the recipe in Section \ref{sec-regeneration},
regeneration times $\tau_i^{(j)}$ corresponding
to the walks $(X_t^{(j)})_{t\geq \bar \beta_m}$ 
and the chain on the environment.
Note that under the measure 
$\bP^{\pi}\otimes 
\bQ_\varepsilon
\otimes
\bQ_\varepsilon\otimes
\Pbar_{\omega^{(1)}, 
\omega^{(1)}, 
{\mathbf \epsilon}^{(1)}, 
{\mathbf \epsilon}^{(2)}}^{\bzero} $, the 
sequences for $j=1$ and $j=2$
are not independent if the walks intersect.
When no confusion occurs, we let $P^{\tau}$ denote the law of this 
sequence,
 noting that the sequence $(\tau_{i+1}^{(j)}-\tau_i^{(j)})_{i\geq 1}$
is i.i.d. under all measures involved in our construction, as well as, 
equation (\ref{eq-190805a}) holds. Fix an integer $k \geq 1$,
for any fixed $2 \leq \gamma' < \gamma$ we then have
\begin{eqnarray}
&      & P^\tau\left(\Big\vert
         \tau^{(j)}_k - k E^\tau [\tau^{(j)}_1] \Big\vert> k/2\right)
         \nonumber \\
&  =   & P^{\tau}\left( \Big\vert \tau_k^{(j)} - k E^{\tau} [\tau_1^{(j)}] 
         \Big\vert^{\gamma'} > \left(k/2\right)^{\gamma'} \right) 
         \nonumber \\
& \leq & \frac{2^{\gamma'} \, 
               E^{\tau}\left[\Big\vert 
               \tau_k^{(j)} - k E^{\tau} [\tau_1^{(j)}] 
               \Big\vert^{\gamma'}\right]}{k^{\gamma'}} \nonumber \\
& \leq & \frac{2^{\gamma'} \, B_{\gamma'} \, 
         E^{\tau}\left[ \Big\vert \tau_1^{(j)} - E^{\tau}[\tau_1^{(j)}] 
         \Big\vert^{\gamma'} \right]}{k^{\gamma'/2}} \,,
         \label{eq-210805bb}
\end{eqnarray}
where the last inequality follows from and application of 
Marcynkiewicz-Zygmund inequality \cite[Pg. 469]{Shi84} 
along with H\"{o}lder inequality,
and $B_{\gamma'} > 0$ is an universal constant which depends on $\gamma'$. 
Recall, $\gamma \theta' / 2 > 1$ and also $\gamma > 6$, thus we can choose
$2 \leq \gamma' < \gamma$ such that $\gamma' \theta' / 2 > 1$. 
Now choosing $k=c_5 m^{\theta'}$, for an appropriate $c_5$, one concludes the
existence of a constant $c_6$ such that
\begin{equation}
\label{eq-210805c}
\sum_{m=1}^\infty
P^\tau\left(
\tau^{(j)}_{\lfloor c_6 m^{\theta'} \rfloor}>
m^{\theta'}/4\right)<\infty\,.
\end{equation}
On the other hand, 
\begin{align*}
&P^\tau\left(\exists i>c_6 m^{\theta'}:\, \tau_{i+1}^{(j)}-
\tau_i^{(j)}>i^\alpha\right)\\
\leq & \sum_{i= \lfloor c_6 m^{\theta'} \rfloor +1}^\infty
P^\tau\left(\tau_{i+1}^{(j)}-
\tau_i^{(j)}>i^\alpha\right)\nonumber \\
= & \sum_{i=\lfloor c_6 m^{\theta'} \rfloor + 1}^\infty
P^\tau\left(\tau_{1}^{(j)} >i^\alpha\right)\\
\leq & \sum_{i= \lfloor c_6 m^{\theta'} \rfloor + 1}^\infty
\frac{C}{i^{\alpha \gamma'}}\,,
\end{align*}
for some deterministic constant $C$, and $\gamma' < \gamma$.
Since our choice of constants
implies that $\theta'(\alpha \gamma -1)>1$, we 
can choose $\gamma' < \gamma$
such that $\theta'(\alpha \gamma' -1)>1$, and so we conclude that
\begin{equation}
\label{eq-210805d}
\sum_{m=1}^\infty
P^\tau\left(\exists i>c_6 m^{\theta'}:\, \tau_{i+1}^{(j)}-
\tau_i^{(j)}>i^\alpha\right)<\infty
\end{equation}
Let $C_m^{(j)}=\{\tau^{(j)}_{ \lfloor c_6 m^{\theta'} \rfloor}
\leq m^{\theta'}/4\} \cap 
\{\forall i>c_6 m^{\theta'}, \, \tau_{i+1}^{(j)}-
\tau_i^{(j)}\leq i^\alpha\}$, and note that
for all measures $P^\tau$ involved,
\begin{equation}
\label{eq-210805e}
\sum_{m=1}^\infty
P^\tau\left(\left(C_m^{(1)}\cap C_m^{(2)}\right)^c\right)
<\infty\,.
\end{equation}
Next, for $j=1,2$, define
$S_t^{(j)}=\widetilde{X}^{(j)}_{\tau_t^{(j)}}$. Recall that 
$\widetilde{X}_s^{(j)}$ is
the orthogonal projection of the vector $X_s^{(j)}$ into the
space $\Lbold^{d-1}$. 
Note that   
\begin{align}
\label{eq-210805f}
&\{
X^{(1)}_{[\bar\beta_m,\infty)}\cap X^{(2)}_{[\bar \beta_m,\infty)}
\neq \emptyset\}\cap A_m\cap C_m^{(1)}
\cap C_m^{(2)} \\
&\quad \quad \quad \subset \{
|S_0^{(1)}-S_0^{(2)}|\geq m^{\theta'},
\exists \ell, k: |S_\ell^{(1)}-S_k^{(2)}|\leq \ell^\alpha+k^\alpha\}
\nonumber 
\end{align}
Let $Z_t$   denote 
a  sum of (centered) $t$ i.i.d. random vectors, 
each distributed according to $S_1^{(1)}$, thus taking values in
$\Lbold^{d-1}$. 
Let $(Z_t')_{t\geq 1}$ denote an independent copy
of $(Z_t)_{t\geq 1}$. Let $P^{\bx,\by}_{Z}$ denote the
law of the sequences $((\bx+Z_t),(\by+Z_t'))_{t\geq 1}$, where
$\bx, \by \in \Lbold^{d-1}$. 
Using that under all the 
measures involved in our computation, the
paths $X^{(1)}$ and $X^{(2)}$ are independent 
until they intersect,  
one concludes that
\begin{align}
\label{eq-240805sleep}
& 
\bP^{\pi}\otimes 
\bP^{\pi}\otimes 
\bQ_\varepsilon
\otimes
\bQ_\varepsilon
\\
&\left( \Pbar_{\omega^{(1)}, 
\omega^{(1)}, 
{\mathbf \epsilon}^{(1)}, 
{\mathbf \epsilon}^{(2)}}^{\bzero} 
\left(\{
X^{(1)}_{[\bar\beta_m,\infty)}\cap X^{(2)}_{[\bar \beta_m,\infty)}
\neq \emptyset\}\cap A_m\cap C_m^{(1)}
\cap C_m^{(2)}\right) \right)\nonumber\\
\leq &
 \mathop{\max_{|\bx|\geq m^{\theta'}}}\limits_{\bx \in \Lbold^{d-1}} 
P_Z^{\bzero,\bx}
\left(\exists \ell,k:\,
|Z_\ell-Z_k'|\leq \ell^\alpha+k^\alpha\right)\,,
\nonumber
\end{align}
with exactly the same estimate when 
$ \Pbar_{\omega^{(1)}, 
\omega^{(2)}, 
{\mathbf \epsilon}^{(1)}, 
{\mathbf \epsilon}^{(2)}}^{\bzero} $ replaces 
$ \Pbar_{\omega^{(1)}, 
\omega^{(1)}, 
{\mathbf \epsilon}^{(1)}, 
{\mathbf \epsilon}^{(2)}}^{\bzero} $.  
But, for some constant $C_\alpha$
 whose value may change from line to line,
\begin{align}
\label{Equ:Large-and-Small-Index-Breaking}
&  \mathop{\max_{|\bx|\geq m^{\theta'}}}\limits_{\bx \in \Lbold^{d-1}}
P_Z^{\bzero,\bx}
\left(\exists \ell, k:\,
|Z_\ell-Z_k'|\leq \ell^\alpha+k^\alpha\right)
\nonumber\\
\leq &
\mathop{\max_{|\bx|\geq m^{\theta'}}}\limits_{\bx \in \Lbold^{d-1}}
P_Z^{\bzero,\bx}
\left(\exists \ell,k:\, |Z_\ell-Z_k'|\leq C_\alpha(\ell+k)^\alpha\right)
\nonumber\\
\leq &
\mathop{\max_{|\bx|\geq m^{\theta'}}}\limits_{\bx \in \Lbold^{d-1}}
P_Z^{\bzero,\bx}
\left(\exists \ell, k: \,(\ell+k)<m^{2\theta''},\,
|Z_\ell-Z_k'|\leq C_\alpha(\ell+k)^\alpha\right) \nonumber \\
&  +
\mathop{\max_{|\bx|\geq m^{\theta'}}}\limits_{\bx \in \Lbold^{d-1}}
P_Z^{\bzero,\bx}
\left(\exists \ell, k: \,(\ell+k)\geq m^{2\theta''},\,
|Z_\ell-Z_k'|\leq C_\alpha(\ell+k)^\alpha\right) \nonumber \\
=: & I_1(m)+ I_2(m)\,.
\end{align}
For $I_1(m)$ 
in (\ref{Equ:Large-and-Small-Index-Breaking}) we observe
that  because $\alpha<1/2$, for any $\gamma'<\gamma$, 
we can find a constant $C_{\gamma'}$ such that for large $m$
\begin{align}
\label{eq-240805a}
& \mathop{\max_{|\bx|\geq m^{\theta'}}}\limits_{\bx \in \Lbold^{d-1}}
P_Z^{\bzero,\bx}
\left(\exists \ell, k: \,(\ell+k)<m^{2\theta''},\,
|Z_\ell-Z_k'|\leq C_\alpha(\ell+k)^\alpha\right)
\nonumber \\
\leq & \, 2 \, P_Z^{\bzero,\bzero}\left(\max_{\ell\leq m^{2\theta''}} |Z_\ell|
\geq \frac{m^{\theta'}}{3} \right) \nonumber \\
\leq & \, \frac{C_{\gamma'}}{m^{(\theta'-\theta'') \gamma'}}\,,
\end{align}
where the first inequality follows from the observation that 
\[
\vert Z_{\ell} - Z_k' \vert \geq \vert \bx \vert - \vert Z_{\ell} \vert
- \vert Z_k' - \bx \vert \,,
\]
and the last inequality is due again to
the Marcynkiewicz-Zygmund inequality coupled with
the Martingale maximal inequality of Burkholder
and Gundy \cite[Pg. 469]{Shi84}.
Therefore, since $\theta''$ was chosen such that
$(\theta'-\theta'')\gamma>1$, it follows that 
\begin{equation}
\label{eq-240805gg}
\sum_m I_1(m)<\infty\,.
\end{equation}

%
Turning to the estimate of the term
$I_2(m)$ in (\ref{Equ:Large-and-Small-Index-Breaking}),
we observe,
with the value of the constants $C_\alpha,C_\alpha'$ possibly
changing from line to line, 
that
\begin{align*}
 & \mathop{\sum_{\ell\geq  k}}\limits_{\left(\ell + k \right) 
\geq m^{2 \theta''}}
\mathop{\max_{|\bx|\geq m^{\theta'}}}\limits_{\bx \in \Lbold^{d-1}}
P_Z^{\bzero,\bx}
\left(|Z_\ell-Z_k'|\leq C_\alpha(\ell+k)^\alpha\right) \\
= & \mathop{\sum_{\ell\geq  k}}\limits_{\left(\ell + k \right) 
\geq m^{2 \theta''}}
\mathop{\max_{|\bx|\geq m^{\theta'}}}\limits_{\bx \in \Lbold^{d-1}}
E_Z^{\bzero, \bx} \left[ P_Z^{\bzero,\bx}
\left(|Z_\ell-Z_k'|\leq C_\alpha(\ell+k)^\alpha \,\Big\vert\, Z_k'\right) 
\right] \\
\leq & 
\mathop{\sum_{\ell\geq  k}}\limits_{\left(\ell + k \right) 
\geq m^{2 \theta''}}
\mathop{\max_{|\bx|\geq m^{\theta'}}}\limits_{\bx \in \Lbold^{d-1}} \, 
C_{\alpha}' \left(\ell + k \right)^{\alpha (d-1)} \,
E_Z^{\bzero, \bx} \left[ 
\max_{\vert \bz - Z_k' \vert \leq C_{\alpha} (\ell + k)^{\alpha}}
P_Z^{\bzero, \bx}\left( Z_{\ell} = \bz \right) \right] \\
\leq & 
\mathop{\sum_{\ell\geq  k}}\limits_{\left(\ell + k \right) 
\geq m^{2 \theta''}}
C_{\alpha}' \left(\ell + k \right)^{\alpha (d-1)} \,
\frac{1}{\ell^{(d-1)/2}} \\
\leq &
\mathop{\sum_{\ell\geq k}}\limits_{\left(\ell + k \right) 
\geq m^{2 \theta''}}
C_{\alpha}' \left(\ell + k \right)^{\alpha (d-1)} \,
\frac{1}{(\ell+k)^{(d-1)/2}} \,,
\end{align*}
where the last but one inequality follows by applying the local limit 
estimate \cite[Theorem 22.1]{BhaRao76}. 
By symmetry, the same estimate holds with 
$k\geq \ell$ replacing 
the constraint 
$\ell\geq k$ in the summation.
Thus we conclude that 
\begin{eqnarray}
I_2\left(m\right) & \leq & 
\mathop{\sum_{\ell, k}}
\limits_{\left(\ell + k \right) \geq m^{2 \theta''}}
\mathop{\max_{|\bx|\geq m^{\theta'}}}\limits_{\bx \in \Lbold^{d-1}}
P_Z^{\bzero,\bx}
\left(|Z_\ell-Z_k'|\leq C_\alpha(\ell+k)^\alpha\right) \nonumber \\
& \leq & 
2 \mathop{\sum_{\ell\geq  k}}\limits_{\left(\ell + k \right)
 \geq m^{2 \theta''}}
C_{\alpha}' \left(\ell + k \right)^{\alpha (d-1)} \,
\frac{1}{\left(\ell + k \right)^{(d-1)/2}} \nonumber \\
& \leq & \sum_{r = \lfloor m^{ 2\theta''} \rfloor + 1}^{\infty} 
\frac{2 C_{\alpha}'}{r^{\sfrac{d-1}{2} - \alpha (d-1) - 1}} 
\leq  
\frac{2 C_{\alpha}'}{m^{\theta''(d-2\alpha d + 2\alpha -5)}}\,.
\label{Equ:Second-Sum-Estimate} 
\end{eqnarray}
In particular, since $\theta''(d-2\alpha d + 2\alpha - 5)>1$, 
one concludes that
\begin{equation}
\label{eq-240805ggg}
\sum_{m} I_2(m)<\infty\,.
\end{equation}
Combining 
(\ref{eq-240805gg}) and (\ref{eq-240805ggg}),
and substituting in (\ref{Equ:Large-and-Small-Index-Breaking})
and then in (\ref{eq-240805sleep}), the conclusion of the lemma follows.

Finally for the case $\bv = \bzero$ we can proceed in two ways. 
We could chose $V^{\bot}$ as an arbitrary subspace
of co-dimension $1$, and simply repeat the argument.
Alternatively, we could proceed exactly the same manner 
except that one works directly with the walkers, rather than their 
projections. In this case we can  replace $(d-1)$ by $d$ in all
the estimates above and the conclusions follows by 
observing that the inequalities in 
(\ref{cond-constants}) holds when $(d-1)$ is
replaced by $d$. The advantage of the second approach is then that
when $\bv=\bzero$, one may replace the constraint on $d$
in (\ref{eq-200805a}) by the weaker constraint involving $d+1$.

This completes the proof of the lemma. $\qed$

\section{Final Remarks and Open Problems}
\label{Sec:Remarks}
\subsection{The ``regeneration'' time}
\label{Subsec:Regeneration}
We point out that our definition of the
``regeneration'' time $\tau_1$ using 
the coin tosses
$\left\{\left(\alpha_t\left(\cdot\right)\right)_{t \geq 1}\right\}$, 
is quite arbitrary, and was tailored to 
Condition (A1).
What we actually need for the argument 
is for every environment chain 
$\left(\omega_t\left(\cdot, \cdot\right)\right)_{t \geq 0}$ a 
sequence of stopping times such that at these times,
the
chain starts from a stationary distribution and the times have 
``good'' tail property. The assumption (A1) is used to ensure
this is possible with our construction. 

\subsection{The annealed CLT}
\label{subsec:remCLTa}
Our argument gives a trade-off between the strength of the random
perturbation and the mixing rate of the environment. We suspect
that such a trade-off is not needed, and that for $d\geq 3$, an annealed
CLT holds true as soon as the environment is mixing enough in time.
Our technique does not seem to resolve this question.

\subsection{The quenched CLT}
\label{subsec:remCLTq}
We have already emphasized that the condition
(\ref{eq-200805a}) is not expected to be optimal. In particular,
we suspect that if a critical dimension for the quenched CLT
exists, it will be $d=1$ (see also the comments 
at the end of \cite{BMP00} hinting that such quenched CLT fails
even for small perturbations of a fixed Markovian environment, and
the numerical simulatons in \cite{BCFP05} that are 
inconclusive).

\section*{Acknowledgment} We thank Amir Dembo for a remark 
that prompted us to modify our original construction of regeneration
times. We are also grateful to the anonymous referee for a careful 
reading of the manuscript and for detecting a mistake in
the original version.


\end{document}